\documentclass{jktr}

\begin{document}

\markboth{V. A. Manathunga}
{The Coefficients of the Jones Polynomial}

\catchline{}{}{}{}{}

\title{THE COEFFICIENTS OF THE JONES POLYNOMIAL}

\author{V. A. Manathunga}

\address{Department of Mathematical Sciences\\
Middle Tennessee State University,\\
 Murfreesboro TN 37032 \\
vmanathunga@mtsu.edu}

\maketitle

\begin{abstract}
It has been known that, the coefficients of the series expansion of the Jones polynomial evaluated at $e^x$  are rational valued Vassiliev invariants . In this article, we calculate minimal multiplying factor, $\lambda$, needed for these rational valued invariants to become integer valued Vassiliev invariants. By doing that we obtain a set of integer-valued Vassiliev invariants.   
\end{abstract}

\keywords{Jones Polynomial, Vassiliev Invariant}

\ccode{Mathematics Subject Classification 2000: 57M25, 57M27}

\section{Introduction}	
 It has been known that coefficients of certain knot polynomials are Vassiliev invariants. For example the $n^{th}$ coefficient of the Conway polynomial is a Vassiliev invariant of order $n$ \cite{dnatan:01}. Vassiliev invariants are an important tool to distinguish knots. It has been conjectured, called \textit{Vassiliev conjecture}, that the set of all rational valued Vassiliev invariants are complete. If this is true, then we can use rational valued Vassiliev invariants to approximate any other rational valued knot invariant. In otherwords if Vassiliev conjecture is true then any link invariant is a pointwise limit of sequence of Vassiliev invariants.  When we expand the Jones polynomial evaluated at $e^x$ using power series representation of $e^x$, the coefficients become rational valued Vassiliev invariants \cite{dnatan:01, goussarov, birman2}. Using this fact, Kofman and Rong \cite{rong} proved that each coefficient of the Jones polynomial of a knot is the limit of a sequence of Vassiliev invariants. Extending this result, 
 Helme-Guizon \cite{helme} approximated the coefficents of HOMFLYPT and Kauffman polynomials using Vassiliev invariants. Turning back to the series expansion of the Jones polynomial, these coefficients are not additive under connected sum, or in other words are not primitive Vassiliev invariants. Since the Jones polynomial, $J_k(t)$, is multiplicative under connected sum, we can take the logarithm of $J_k(e^x)$ and expand it as a formal power series. However this turns integer coefficients of Jones polynomial, $J_k(t)$, into rational valued ones.  Whenever we have rational valued Vassiliev invariants we can construct an integer valued invariant by multiplying by a suitable ``large'' integer. However if we do that resulting invariants tend to have common factors. Here we calculate minimal multiplying factor, $\lambda$, needed for the coefficient of the Jones polynomial evaluated at $e^x$ and expanded using power series representation of $e^x$ to become integer valued.

\section{Background}

\begin{definition}\label{defn20}
	The Jones polynomial $J_K(t)$ is defined as the polynomial satisfying the following \em{skein relation}
	\begin{gather}
		t^{-1}J_K\Bigl( \text{\raisebox{-6pt}[0pt][0pt]{\includegraphics{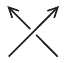}}} \Bigr)-t J_K\Bigl(\text{\raisebox{-6pt}[0pt][0pt]{\includegraphics{Stoimenow-Pic1}}}  \Bigr) = \big(t^{1/2}-t^{-1/2}\big)J_K\Bigl(\text{\raisebox{-5pt}[0pt][0pt]{\includegraphics{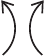}}} \Bigr).\\
		J_K(unknot)=1
	\end{gather}
	
\end{definition}

Let the Jones polynomial of a knot $K$ be $J_K(t)$ where ,
\begin{align*}
	J_k(t)&=\Sigma_{i=-m}^{n}c_it^i \\
	&=c_{-m}t^{-m}+c_{-m+1}t^{-m+1}+...+c_{-1}t^{-1}+c_{0}+c_1t^1+...+c_nt^n \text{\quad where } c_i\in \mathbb{Z}
\end{align*}

Let $t=e^x$ then $J_K(e^x)=\Sigma_{k=0}^{\infty}a_kx^k$ is the Jones polynomial of knot $K$ expanded as a power series in $x$. It is known that $a_0=1,a_1=0$ and for $k>1$, $a_k=\frac{1}{k!}\Sigma_{i=-m,i\ne 0}^nc_ii^k$ is a rational valued Vassiliev invariant. We can convert $a_k$ to be integer valued by multiplying by a suitable constant $\lambda_k$. For example, we can choose $\lambda_k=k!$ and then $\lambda_ka_k$ is an integer valued Vassiliev invariant. In this article we calculate minimal $\lambda_k$ such that $\lambda_ka_k$ is an integer valued Vassiliev invariant. First we define what we meant by ``minimal'' $\lambda_k$. By ``minimal''$\lambda_k$, we mean the unique positive constant $\lambda_k$ which makes the values of $\lambda_ka_k$ integers and makes it so that not all values of $\lambda_ka_k$ have a common factor. The formal definition as follows: 
\begin{definition}
	Let $\K$ denote the set of equivalence classes of knots. Let,
	\begin{align*}
		q_k &=\lcm \{q_i | a_k(K_i)=\frac{p_i}{q_i}, \gcd(p_i,q_i)=1, K_i\in \K  \}\\
		p_k &=\gcd \{p_i | a_k(K_i)=\frac{p_i}{q_i}, \gcd(p_i,q_i)=1, K_i\in \K  \} \text{\quad then,}\\
		\lambda_k &=\frac{q_k}{p_k}
	\end{align*}
\end{definition}
The existence of $q_k$ and $p_k$ over infinite set $\K$ may come in to the question at this point. However it is due to the fact that $a_k$ are Vassiliev invariants and the value of Vassiliev invariant over any given knot is a integer linear combination of its value on only finite number of singular knots.

It is clear from the definition $\lambda_ka_k\in \Z$, $\gcd(p_k,q_k)=1$ and $k!$ is an upper bound of $\lambda_k$. One way to calculate $\lambda_k$ is actuality tables. Suppose $v_n$ is a Vassiliev invariant of order $\le n$. To construct actuality tables one must choose a singular knot for every basis chord diagram of order $\le n$. Then one must calculate the value of a given invariant over all those basis knot diagrams. This can be done easily if $n$ is relatively small. For example, from actuality tables, we can see that $\lambda_2=\frac{1}{3}$ and $\lambda_3=\frac{1}{6}$. However this method is not feasible when $k\ge 4$. So, from this point onward we are interested about $a_k,$ and $\lambda_k$, only when $k\ge 4$. The following lemma gives more insight to the above definition.

\begin{lemma}\label{lemma1}
	Suppose there exists $\lambda'_k$ and two knots $K_1, K_2\in \K$ such that,
	\begin{enumerate}
		\item $\lambda'_ka_k\in \mathbb{Z}$
		\item $\gcd(\lambda'_ka_k(K_1),\lambda'_ka_k(K_2))=1$ 
	\end{enumerate}
	Then $\lambda'_k$ is minimal.
\end{lemma}

\begin{proof}
	Let $\lambda_k=\frac{q_k}{p_k}$ and $\lambda'_k=\frac{q'_k}{p'_k}$, where $\gcd(p_k,q_k)=1,\gcd(p'_k,q'_k)=1$. Then from the definition of $\lambda_k$, it is clear that $q_k|q'_k$ and $p'_k|p_k$. We can easily prove that, $\lambda'_k=r\lambda_k$ for some $r\in \Z$. Since $\lambda_ka_k(K_1),\lambda_ka_k(K_2)\in \mathbb{Z}$, we conclude $\gcd(\lambda'_ka_k(K_1),\lambda'_ka_k(K_2))=r=1$. Hence $\lambda'_k=\lambda_k$.
\end{proof}
The following lemma lists some known results about coefficients of Jones polynomial. We use this lemma extensively in subsequent theorems.
\begin{lemma}\label{lemma2}
	Let Jones polynomial of knot $K$ be $J_K(t)=\Sigma_{i=-m}^{n}c_it^i$. Then,
	\begin{enumerate}
		\item $J_K(1)=\Sigma_{i=-m}^{n}c_i=1$
		\item $J_K(-1)=\Sigma_{i=-m}^{n}(-1)^ic_i\equiv 1 \text{ or } 5 \mod 8$
		\item $J_K(e^{2\pi i/ 3})=1$
		\item $J_K(\sqrt{-1})=(-1)^{c_2(K)}$ where $c_2(K)$ is the second coefficient of the Conway polynomial of knot $K$.
		\item $J_K(1)-J_K(-1)=\Sigma_{i=-m}^{n}(c_i-(-1)^ic_i)\equiv 0\mod 4$ which implies \\ $\Sigma_{i=-m,2\nmid i}^{n}c_i\equiv 0\mod 2$ 
		\item $\Sigma_{i=-m,3|i }^{n}c_i=\cdot\cdot\cdot+c_{-6}+c_{-3}+c_0+c_3+c_6+\cdot\cdot\cdot=1$
		\item $\Sigma_{i=-m,3\nmid i }^{n}(-1)^{(i-1)\mod 3}c_i=\cdot\cdot\cdot-c_{-4}+c_{-2}-c_{-1}+c_1-c_2+c_4-\cdot\cdot\cdot=0$
		\item $\Sigma_{i=-m,3\nmid i }^{n}c_i=\cdot\cdot\cdot+c_{-5}+c_{-4}+c_{-2}+c_{-1}+c_1+c_2+c_4+c_5+\cdot\cdot\cdot=0$
		\item $\Sigma_{i=-m,2\nmid i }^{n}(-1)^{(i-1)/2}c_i=\cdot\cdot\cdot-c_{-5}+c_{-3}-c_{-1}+c_1-c_3+c_5\cdot\cdot\cdot=0$
	\end{enumerate}
\end{lemma}

\begin{proof}
	We only prove $(6),(7),(8)$ and $(9)$. Others are standard results. Now from (3), we have, 
	$J_K(e^{(2\pi i/ 3)})=\Sigma_{i=-m,3|i }^{n}c_i - \frac{1}{2}\Sigma_{i=-m,3\nmid i }^{n}c_i +\frac{i\sqrt{3}}{2} \Sigma_{i=-m,3\nmid i }^{n}(-1)^{(i-1)\mod 3}c_i=1$.
	This imply, $\Sigma_{i=-m,3|i }^{n}c_i - \frac{1}{2}\Sigma_{i=-m,3\nmid i }^{n}c_i =1 $ and $ \Sigma_{i=-m,3\nmid i }^{n}(-1)^{(i-1)\mod 3}c_i =0$. Using (1), we can rewrite, $\Sigma_{i=-m,3|i }^{n}c_i - \frac{1}{2}\Sigma_{i=-m,3\nmid i }^{n}c_i=\Sigma_{i=-m,3|i }^{n}c_i - \frac{1}{2}(1-\Sigma_{i=-m,3|i }^{n}c_i)=1$. Thus we have $\Sigma_{i=-m,3|i }^{n}c_i=1$. Substituting this result in (1), we get $\Sigma_{i=-m,3\nmid i }^{n}c_i=0$. Similarly, $J_K(\sqrt{-1})=$ $\sqrt{-1}(\Sigma_{i=-m,2\nmid i }^{n}(-1)^{(i-1)/2}c_i)$ + $(c_0+\Sigma_{i=-m,2|i,i\ne 0 }^{n}(-1)^{(i-1)/2}c_i)$ $=(-1)^{c_2(K)}$. This implies the result of (9). 
\end{proof}
Now consider the $k^{th}$ coefficient, $a_k$ of the Jones polynomial evaluated at $e^x$ and expanded using power series representation of $e^x$ where $J_K(t)=\Sigma_{i=-m}^{n}c_it^i$. Then it is clear that when $k>0$, $a_k=\frac{1}{k!}\Sigma_{i=-m,i\ne 0}^nc_ii^k$. The following lemma gives some basic properties of $a_{k}$. 

\begin{lemma}{\ \\}\label{lemma3}
	\begin{enumerate}
		\item $(2k)!a_{2k}\equiv 0 \mod 2$, where $k\ge 1$.
		\item Let $n$ be any non negative integer and $k\ge 2$, If $k\equiv 0\mod 3^{n}$ then $(k)!a_{k}\equiv 0 \mod 3^{n+1}$.
		\item $(2k+1)!a_{2k+1}\equiv 0 \mod 4$, where $k\ge 1$.
	\end{enumerate}
\end{lemma}

\begin{proof}{\ \newline}
	\begin{enumerate}
		\item First observe that, $(2k)!a_{2k}  \mod{2}$ $\equiv$ $\sum_{i=-m, 2\nmid i}^n c_i \mod{2}$. Now from lemma \ref{lemma2}(5) we conclude that $(2k)!a_{2k}\equiv 0$ mod $2$.
		
		\item If $k\equiv 0\mod 3^n$ then, 
		\[
		i^{2k}\equiv
		\begin{cases}
		0\mod 3^{n+1}, & \text{if } i\equiv 0 \mod{3}\\
		1\mod 3^{n+1}, & \text{otherwise}\\
		\end{cases}
		\]
		Hence, $(2k)!a_{2k}$ $\equiv$ $\sum_{i=-m, 3\nmid i}^n c_i\mod{3^{n+1}}$. Now from lemma \ref{lemma2}(8), we conclude that $(2k)!a_{2k}\equiv 0\mod 3^{n+1}$.\\
		Similarly,
		\[
		i^{2k+1} \equiv
		\begin{cases}
		1\mod 3^{n+1}, & \text{if } i\equiv 1 \mod{3}\\
		(-1) \mod 3^{n+1}, & \text{if } i\equiv (-1) \mod{3}\\
		0 \mod 3^{n+1}, & \text{otherwise}\\
		\end{cases}
		\]
		Thus, $(2k+1)!a_{2k+1}$ $\equiv$ $\Sigma_{i=-m,3\nmid i }^{n}(-1)^{(i-1)\mod 3}c_i \mod 3^{n+1}$. Therefore, from lemma \ref{lemma2}(7) we conclude that, $(2k+1)!a_{2k+1}$ $\equiv$ $0\mod 3^{n+1}$. 
		\item  Since $2k+1\ge 3$, 
		\[
		i^{2k+1} \equiv
		\begin{cases}
		1\mod 4, & \text{if } i\equiv 1 \mod{4}\\
		(-1) \mod 4, & \text{if } i\equiv (-1) \mod{4}\\
		0 \mod 4, & \text{otherwise}\\
		\end{cases}
		\]
		Therefore, $(2k+1)!a_{2k+1}$ $\equiv$ $\Sigma_{i=-m,2\nmid i }^{n}(-1)^{(i-1)/2}c_i \mod 4$. Thus from lemma \ref{lemma2}(9) we conclude that $(2k+1)!a_{2k+1}\equiv 0 \mod 4$.
		
	\end{enumerate}
\end{proof} 
\section{Main Theorem }
Now we are ready to prove the main result of this article. In the following theorem we give a formula for $\lambda_{2k}$ when $k\ge 1$. 
\begin{theorem}\label{thm1}
	Let $J_K(e^x)=\sum_{i=0}^{\infty}a_kx^k$ be the Jones polynomial of knot $K$ expanded as a power series in $x$. When $k\ge 1$, let $\lambda_{2k}=\frac{(2k)!}{2.3^{n+1}}$ where $3^n|k$ but $3^{n+1}\nmid k$ and $n$ is a non negative integer. Then,
	\begin{enumerate}
		\item $v_{2k}=\lambda_{2k}a_{2k}(K)$ is an integer valued Vassiliev invariant .
		\item $\lambda_{2k}$ is minimal.
	\end{enumerate}
\end{theorem}

\begin{proof}
	From lemma \ref{lemma3}, we have $2.3^{n+1}|(2k)!a_{2k}$. Therefore, $v_{2k}=\lambda_{2k}a_{2k}\in \mathbb{Z}$. Next we prove $\lambda_{2k}$ is minimal.  We claim $gcd(v_{2k}(3_1),v_{2k}(4_1))=1$ for any given $k$. First observe that $J_{3_1}(t)=-t^{-4}+t^{-3}+t^{-1}$ and $J_{4_1}(t)=t^{-2}-t^{-1}+1-t+t^2$. Let $ 3^n|k,3^{n+1}\nmid k$ then,
	\begin{align*}
		v_{2k}(3_1) &=\left[\frac{(2k)!}{2.3^{n+1}}\right ]\left[ \frac{-(-4)^{2k}+(-3)^{2k}+(-1)^{2k}}{(2k)!}\right] \\ 
		&= \frac{-16^k+9^k+1}{2.3^{n+1}}
	\end{align*}
	Similarly,
	\begin{align*}
		v_{2k}(4_1) &=\left[\frac{(2k)!}{2.3^{n+1}}\right ]\left[ \frac{(-2)^{2k}-(-1)^{2k}-(1)^{2k}+(2)^{2k}}{(2k)!}\right] \\ 
		&= \frac{2(4^k-1)}{2.3^{n+1}}
	\end{align*}
	\\
	We prove,$\gcd(-16^k+9^k+1,2(4^k-1))=2.3^{n+1}$. From above we have $2.3^{n+1}|(2k)!a_{2k}$. Therefore $2.3^{n+1}|\gcd(-16^k+9^k+1,2(4^k-1))$. Now observe that, $2^2\nmid 2(4^k-1)$ thus $2^2\nmid \gcd (-16^k+9^k+1,2(4^k-1))$. Next observe that if,  $3^{n+2}|\gcd (-16^k+9^k+1,2(4^k-1))$ then $3^{n+2}|2(4^k-1)$. This imply $3^{n+2}|4^{k}-1$. But $k=3^n.q$ where $3\nmid q$. Thus we have $3^{n+2}|(3+1)^{3^n.q}-1$. This give us $3^{n+2}|3^{n+1}$, hence a contradiction. 
	
	Let $P$ be a prime such that $P\ne 2,3$ and $P |\gcd(-16^k+9^k+1,2(4^k-1))$. This imply, $P|4^k-1$ and $P|-16^k+9^k+1$. Now observe that, $-16^k+1=-[4^k-1]-[4^k(4^k-1)]$. Therefore, $P|-16^k+1$. Thus $P|9^k$. This give contradiction, since $P\ne 2,3$. Hence we conclude, $\gcd(-16^k+9^k+1,2(4^k-1))=2.3^{n+1}$. Thus $\gcd(v_{2k}(3_1),v_{2k}(4_1))=1$ for any given $k$. So, from Lemma \ref{lemma1}, $\lambda_{2k}$ is minimal for any $k$.
	
\end{proof}
By induction, we can show easily $\lambda_{2k}\in \Z$ when $k\ge 2$. Thus when $k\ge 2$,  
\begin{align*}
	\lcm \{q_i | a_{2k}(K_i)=\frac{p_i}{q_i}, \gcd(p_i,q_i)=1, K_i\in \K  \} &=\frac{(2k)!}{2.3^{n+1}}\text{\quad where $3^n|k$, $3^{n+1}\nmid k$ }\\
	\gcd \{p_i | a_{2k}(K_i)=\frac{p_i}{q_i}, \gcd(p_i,q_i)=1, K_i\in \K  \} &=1
\end{align*}
This is nice fact to know given that family of knots is infinite.\\

Now we should focus on finding a general formula for $\lambda_{2k+1}$. However it seems that this is not as easy as the $\lambda_{2k}$ case. First observe that from Lemma \ref{lemma3}, $\lambda_{2k+1}\le \frac{(2k+1)!}{12}$. Computational results shows that this is indeed the least upper bound for some $k$. For example, $\lambda_{35}=35!/12$, $\lambda_{95}=95!/12$, $\lambda_{119}=119!/12$, $\lambda_{143}=143!/12$. This is because, when $k=35,95,119,143$ we can show $\lcm \{\,q_k(3_1),q_k(5_1),q_k(9_1) \}\,= \frac{k!}{12}$ and $\gcd \{\,p_k(3_1),p_k(5_1),p_k(9_1) \}\,= 1$. The seeming randomness of when $\lambda_{2k+1}=(2k+1)!/12$ to hampers giving a general formula for $\lambda_{2k+1}$.\\

A careful reader may have noticed that the above invariants, $v_k=\lambda_k a_k$, are integer valued but not primitive ones. Let $J_K(e^x)=1+\Sigma_{i=1}^{\infty}a_ix^i$. Then $\log(J_k(e^x))=\Sigma_{i=2}^{\infty}w_ix^i $ where,
\[
|w_{i}|=\sum_{\Sigma k_jn_j=i} \Bigg(\frac{1}{n_1+\cdot\cdot\cdot+n_m}\Bigg){n_1+n_2+\cdot\cdot\cdot+n_m\choose n_1,n_2,\cdot\cdot\cdot,n_m} a_{k_1}^{n_1} \cdot\cdot\cdot a_{k_m}^{n_m}
\]
We conjecture that $\lambda_{2k}w_{2k}$ is a primitive integer valued Vassiliev invariant for all $k>1$ and $\lambda_{2k}$ is ``minimal''. Computational results shows that $\gcd(\lambda_{2k}w_{2k}(3_1),\lambda_{2k}w_{2k}(4_1),\lambda_{2k}w_{2k}(5_1))=1$ for $1\le k\le 300 $ and we conjecture $\gcd(\lambda_{2k}w_{2k}(3_1),\lambda_{2k}w_{2k}(4_1),\lambda_{2k}w_{2k}(5_1))=1$ for all $k\in\mathbb{Z}$. 

Let $\mathbb{Z}_{(2)}$ denote the localization of $\mathbb{Z}$ at $(2)$, where $(2)$ is the prime ideal generate by $2\in \mathbb{Z}$. Then we have following interesting lemma.
\begin{lemma}\label{lemma4} When $k\ge 1$,
	\begin{enumerate}
		\item $2^{k-1}k!a_{2k}\in \mathbb{Z}_{(2)}$
		\item $2^{k-2}k!a_{2k+1}\in \mathbb{Z}_{(2)}$
	\end{enumerate}
\end{lemma}
\begin{proof}
	Using mathematical induction, we can show that $\frac{2^kk!}{(2k)!}\in \Z_{(2)}$. Now $2^{k-1}k!a_{2k}=\frac{2^{k-1}k!}{(2k)!}(2k)!a_{2k}$. Since we proved $(2k)!a_{2k}\equiv 0 \mod{2}$ in Lemma \ref{lemma3}, the first part of the lemma follows. Similarly $2^{k-2}k!a_{2k+1}=\frac{2^{k-2}k!}{(2k+1)(2k)!}(2k+1)!a_{2k+1}$ and since we proved $(2k+1)!a_{2k+1}\equiv 0\mod{4}$ in Lemma \ref{lemma3}, second part of the lemma follows.
\end{proof}
\newpage
The following lemma shows an interesting characteristic of $\lambda_ka_k \mod 2$.
\begin{lemma}\label{lemma5} Let $v_k=\lambda_ka_k$. Then $v_{2k}\equiv v_{2k+2}\mod 2$ for all $k\ge 1$.
\end{lemma} 
\begin{proof}
	First observe that $v_{2k}=\frac{(2k)!}{2.3^{n+1}}\Sigma_{i=-m,i\ne 0}^{n'}\frac{c_ii^{2k}}{(2k)!}\in \mathbb{Z}$. Thus $v_{2k}3^{n+1}=\frac{1}{2}\Sigma_{i=-m,i\ne 0}^{n'}c_ii^{2k}$ is integer valued and we conclude, $v_{2k}\equiv \frac{1}{2}\Sigma_{i=-m,i\ne 0}^{n'}c_ii^{2k} \mod 2 $.  Now $v_{2k+2}-v_{2k}\equiv \Sigma_{i=-m,i\ne 0}^{n'}c_i(\frac{i^{2k+2}-i^{2k}}{2}) \mod 2 $. It is easy to prove that $i^{2k+2}-i^{2k}\equiv 0\mod 4$. Hence $v_{2k+2}
	\equiv v_{2k}\mod 2$. 
\end{proof}

If the following conjecture is true, then we can extend the lemma \ref{lemma5} to odd case using a similar type proof.
\begin{conjecture}
	When $k\ge 1$, $(2k+1)!a_{2k+1}\not\equiv 0 \mod{8}$
\end{conjecture}

Now, we give $\lambda_{2k}$ values and conjectured $\lambda_{2k+1}$ values. 
\begin{table}
	\centering
	\caption{Values of $\lambda_{2k}$ and conjectured $\lambda_{2k+1}$ up to order $15$}
	\begin{tabular}{|l|l|}
		\hline $\lambda_2=\frac{1}{3}=\frac{2!}{2\times 3}$ & $\lambda_3=\frac{1}{6}=\frac{3!}{2^2\times 3^2}$ \\ 
		\hline $\lambda_4=4=\frac{4!}{2\times 3}$ & $\lambda_5=2=\frac{5!}{2^2\times 3\times 5}$  \\ 
		\hline $\lambda_6=40=\frac{6!}{2\times 3^2}$  & $\lambda_7=60=\frac{7!}{2^2\times 3\times 7}$ \\ 
		\hline $\lambda_8=6720=\frac{8!}{2\times 3}$ & $\lambda_9=672=\frac{9!}{2^2\times 3^3\times 5}$ \\ 
		\hline $\lambda_{10}=604800=\frac{10!}{2\times 3}$ & $\lambda_{11}=302400=\frac{11!}{2^2\times 3\times 11}$ \\ 
		\hline $\lambda_{12}=26611200=\frac{12!}{2\times 3^2}$ & $\lambda_{13}=1140480=\frac{13!}{2^2\times 3\times 5\times 7\times 13}$ \\ 
		\hline $\lambda_{14}=14529715200=\frac{14!}{2\times 3}$ & $\lambda_{15}=36324288000=\frac{15!}{2^2\times 3^2}$  \\ \hline 
	\end{tabular}
	
\end{table}

These integer valued Vassiliev invariants seemed to satisfies multiple congruence relations. 

\begin{proposition}{\ \\}
	Let Jones polynomial of knot $K$ be $J_K(t)=\Sigma_{i=-m}^{n}c_it^i$. Assume, $v_7, v_9, v_{11}$ and $v_{15}$ are integer valued and $\lambda_7, \lambda_9, \lambda_{11}$ and $\lambda_{15}$ as in the above table. Then,
	\begin{enumerate}
		\item $v_4\equiv v_8\mod{5}$
		\item $v_4\equiv v_{10}\mod{7}$
		\item $v_6\equiv v_{12}\mod{7}$
		\item $v_3\equiv v_{15}\mod{13}$
		\item $v_3\equiv v_9\mod{7}$
		\item $v_7\equiv v_{11}\mod{24}$
	\end{enumerate}
	
\end{proposition}

\begin{proof}
	First observe that,
	\[
	\frac{\lambda_4}{4!}=\frac{\lambda_8}{8!}, \frac{\lambda_4}{4!}=\frac{\lambda_{10}}{10!}, \frac{\lambda_6}{6!}=\frac{\lambda_{12}}{12!}, \frac{\lambda_3}{3!}=\frac{\lambda_{15}}{15!}
	\]
	Now $v_4-v_8=\frac{1}{6}\sum_{i=-m,i\ne 0}^nc_i(i^4-i^8)$. Since we proved $v_{2k}$ is an integer, 
	\[6(v_4-v_8)\mod{5}\equiv (v_4-v_8)\mod{5}\equiv \sum_{i=-m,i\ne 0}^nc_ii^4(1-i^4)\mod{5}.\] But $i^4(1-i^4)\equiv 0\mod{5}$. Hence the result follows. Similar proofs can be given for (2) and(3).
	
	For (4), \[v_{15}-v_3=\frac{1}{2^23^2}\sum_{i=-m,i\ne 0}^nc_i(i^{15}-i^3)=\frac{1}{2^23^2}\sum_{i=-m,i\ne 0}^nc_ii^3(i^{12}-1).\] But $i^{12}-1\equiv 0 \mod{13}, i^3(i^{12}-1)\equiv 0 \mod{4}, i^3(i^{12}-1)\equiv 0 \mod{9}$. Hence result follows.
	
	Assume, $v_9, v_7$ and $v_{11}$ are integer valued invariants. From actuality tables we know that $v_3$ is an integer valued.
	
	For(5) First observe that $2^23^2v_3\equiv v_3\mod{7}$ and $2^23^35 v_{9}\equiv v_9\mod{7}$. Now $(v_9-v_3)\equiv\sum_{i=-m,i\ne 0}^nc_i (i^9-i^3)\mod{7}$. But it is easy to prove that $i^9-i^3\equiv 0\mod{7}$. Hence result follows. 
	
	For(6), Similar proof can be given.
\end{proof}

 Before we end this article, we state following two lemmas showing two characteristics of $\lambda_{2k}$, which we calculated above .

\begin{lemma}
	$\lambda_{2i}$ is divisible by $\lambda_j$ for all $j\le 2i$
\end{lemma}

\begin{proof}
	Suppose $j=2k$ and assume, $\lambda_{2i}=\frac{(2i)!}{2.3^{n+1}}, 3^n|i, 3^{n+1}\nmid i$ and $\lambda_{2k}=\frac{(2i)!}{2.3^{m+1}}, 3^m|k, 3^{m+1}\nmid k$ where $m,n\in \N0$ .
	\begin{align*}
		\frac{\lambda_{2i}}{\lambda_{2k}}&=\frac{(2i)!3^{m+1}}{(2k)!3^{n+1}}\\
		&=\frac{(2i).q}{3^{n-m}} \text{ for some } q\in \mathbb{Z}
	\end{align*}
	If $n\le m$, then it is clear that,$\lambda_{2k}|\lambda_{2i}$. Suppose $n>m$. Then,
	\[
	\frac{\lambda_{2i}}{\lambda_{2k}}=\frac{(2.3^n.t_1).q}{3^{n-m}}\qquad  \text{ where } i=3^n.t_1
	\]
	Thus, it is clear that $\lambda_{2k}|\lambda_{2i}$. 
	Now suppose, $j={2k+1}$. Then, $\lambda_{2k+1}=\frac{(2k+1)!}{2.3^{m+1}.r}$ where $3^m|2k+1,3^{m+1}\nmid 2k+1$, $m\in \N0$and $r\in \mathbb{Z}$ . 
	\begin{align*}
		\frac{\lambda_{2i}}{\lambda_{2k+1}}&=\frac{(2i)!2.3^{m+1}r}{(2k+1)!3^{n+1}} \\
		&=\frac{(2i)q_2}{3^{n-m}} \text{ for some } q_2\in \Z
	\end{align*} 
	Now using same argument as above, we can conclude $\lambda_{2k+1}|\lambda_{2i}$. Thus $\lambda_j|\lambda_{2i}$ for all $j\le 2i$
\end{proof}

\begin{lemma}
	$\prod_{\Sigma_jk_j=2i} \lambda_{k_j}|\lambda_{2i}$ where $2\le k_j\le 2i-2$ .
\end{lemma}
\begin{proof}
	From theorem \ref{thm1}, we have $\lambda_{2i}=\frac{(2i)!}{2.3^{n+1}}$ where $3^n|i$ but $3^{n+1}\nmid i$ and $n\in \N0$. Now
	\begin{align*}
		\frac{\lambda_{2i}}{\prod_{\Sigma_jk_j=2i}\lambda_{k_j}}&=\Bigg(\frac{(2i)!}{\prod_{\Sigma_jk_j=2i}(k_j)!} \Bigg)\Bigg( \frac{q_1}{2.3^{n+1}}\Bigg) \text{for some } q_1\in\mathbb{N}
	\end{align*}
	Observe that $\frac{(2i)!}{\prod_{\Sigma_jk_j=2i}(k_j)!}$ denote multinomial coefficients, which is a generalization of binomial coefficients and thus an integer.
	Since $\lambda_{2i-2}=\frac{(2i-2)!}{2.3^{m+1}}$ for some $m\in \N0$ we can see $q_1=2.3^{m+1}q_2$ for some $q_2\in \mathbb{N},m\in \N0$.  Thus we have to consider only the term $3^{n+1}$ in the denominator of the $\Bigg(\frac{q_1}{2.3^{n+1}}\Bigg)$ in order to claim $\frac{\lambda_{2i}}{\prod_{\Sigma_jk_j=2i}\lambda_{k_j}}$ is an integer. If $n=0$, then since the minimum $m$ can take is 0, we are done. Otherwise we know that $3|i$. Suppose $\frac{i}{3}\ne 2$. Then consider the $\lambda_2\lambda_{i/3}=\Big(\frac{(2)!}{2\times 3}\Big)\Big( \frac{(i/3)!}{2^2\times 3^n\times q_3}\Big)$ for some $q_3\in\mathbb{N}$. Thus we can see that by dividing $\lambda_{2i}$ by $\lambda_{2}\lambda_{i/3}$ we have $q_1=(2\times 3)(2^2\times 3^n\times q_4)$ for some $q_4\in \mathbb{N}$. Therefore we cancel, $3^{n+1}$ term in the denominator of $\Bigg(\frac{q_1}{2.3^{n+1}}\Bigg)$. Now suppose $i/3=2$. This implies $2i=12$. In that case $\lambda_{12}=\frac{12!}{2.3^2}$. Now using $\lambda_{6}=\frac{6!}{2\times 3^2}$ we can cancel the $2.3^2$ term in the denominator of the $\lambda_{12}$. Hence result follows  
\end{proof}

\section{Acknowledgment}

The author wish to thank Jim Conant for helpful discussions.

\end{document}